\documentclass[11pt]{article}

\usepackage{mathrsfs}
\usepackage{amsfonts}
\usepackage{}

\usepackage[leqno]{amsmath}
\usepackage{graphicx}
\usepackage{latexsym}
\usepackage{amsmath,amsfonts,amssymb,amsthm,mathrsfs,euscript,makeidx,color}
\usepackage{enumerate}

\allowdisplaybreaks[2]

\def\sqr#1#2{{\vcenter{\vbox{\hrule height.#2pt
				\hbox{\vrule width.#2pt height#1pt \kern#1pt \vrule width.#2pt}
				\hrule height.#2pt}}}}
\def\5n{\negthinspace \negthinspace \negthinspace \negthinspace \negthinspace }
\def\4n{\negthinspace \negthinspace \negthinspace \negthinspace }
\def\3n{\negthinspace \negthinspace \negthinspace }
\def\2n{\negthinspace \negthinspace }
\def\1n{\negthinspace }

\def\dbE{\mathbb{E}}

\def\dbP{\mathbb{P}}

\def\dbS{\mathbb{S}}

\def\cF{{\cal F}}

\def\ds{\displaystyle}

\def\ns{\noalign{\ss}}

\def\q{\quad}
\def\qq{\qquad}

\def\({\Big (}
\def\){\Big )}
\def\[{\Big[}
\def\]{\Big]}

\def\a{\alpha}
\def\b{\beta}

\def\l{\lambda}

\def\t{\tau}

\def\th{\theta}

\def\cd{\cdot}

\def\bde{\begin{definition}\label}
	\def\ede{\end{definition}}
	\def\bel{\begin{equation}\label}
		\def\ee{\end{equation}}
	\def\bt{\begin{theorem}\label}
		\def\et{\end{theorem}}
	\def\bc{\begin{corollary}\label}
		\def\ec{\end{corollary}}
	\def\bl{\begin{lemma}\label}
		\def\el{\end{lemma}}
	\def\bp{\begin{proposition}\label}
		\def\ep{\end{proposition}}
	\def\bas{\begin{assumption}\label}
		\def\eas{\end{assumption}}
	\def\br{\begin{remark}\label}
		\def\er{\end{remark}}
	\def\bex{\begin{example}\label}
		\def\ex{\end{example}}
	\def\ba{\begin{array}}
		\def\ea{\end{array}}
	\def\ben{\begin{enumerate}}
		\def\een{\end{enumerate}}
	
	\def\square#1{\vbox{\hrule\hbox{\vrule height#1%
				\kern#1\vrule}\hrule}}
	\def\rectangle#1#2{\vbox{\hrule\hbox{\vrule height#1%
				\kern#2\vrule}\hrule}}

	\font\tenbb=msbm10 \font\sevenbb=msbm7 \font\fivebb=msbm5
	
	\newfam\bbfam
	\scriptscriptfont\bbfam=\fivebb \textfont\bbfam=\tenbb
	\scriptfont\bbfam=\sevenbb

	\newtheorem{theorem}{\indent Theorem}[section]
	\newtheorem{definition}[theorem]{\indent Definition}
	\newtheorem{proposition}[theorem]{\indent Proposition}
	\newtheorem{corollary}[theorem]{\indent Corollary}
	\newtheorem{lemma}[theorem]{\indent Lemma}
	\newtheorem{remark}[theorem]{\indent Remark}
	\newtheorem{example}[theorem]{\indent Example}
	
	\newtheorem{assumption}[theorem]{\indent Assumption}
	
	\makeatletter
	
	\@addtoreset{equation}{section}
	\makeatother
	
	\def\bea{\begin{equation*}}
		\def\eea{\end{equation*}}
	
	\def\bel{\begin{equation}\label}
		\def\eel{\end{equation}}
	
	\def\ba{\begin{array}}
		\def\ea{\end{array}}
	\newcommand{\ad}{&\!\!\!\displaystyle}

	\def\({\Big (}
	\def\){\Big )}
	\def\[{\Big[}
	\def\]{\Big]}
	\def\q{\quad}
	\def\qq{\qquad}

	\def\ds{\displaystyle}

	\def\ns{\noalign{\smallskip}}

\def\PP{\mathsf P}

\def\EE{\mathsf E}

\begin{document}


\title{Minimax Sequential Testing for Poisson Processes}

\author{
Hongwei Mei\thanks{ Department of Mathematics and Statistics, Texas Tech University, Lubbock, TX 79409, USA; email: {\tt hongwei.mei@ttu.edu.} This author is partially supported by Simons Foundation's Travel Support for Mathematicians Program (No. 00002835).}  ~~~ and ~~~ Rui Wang\thanks{ Department of Mathematics and Statistics, Texas Tech University, Lubbock, TX 79409, USA; email: {\tt rui-math.Wang@ttu.edu.} }
}

\maketitle

\begin{abstract}
Suppose we observe a Poisson process in real time for which the intensity may take on two possible values $\lambda_0$ and $\lambda_1$. Suppose further that the priori probability of the true intensity is not given.  We solve a minimax version of  Bayesian problem of sequential testing of two simple hypotheses to minimize a linear combination of the probability of wrong detection and the expected waiting time in the worst scenario of all possible priori distributions. An equivalent characterization for the least favorable distributions is derived and a sufficient condition for the existence is concluded.
\end{abstract}

\paragraph{MSCcodes}Primary: 60G40, 62C20. Secondary: 60H30.

\section{Introduction}
	Suppose that we observe a continuous-time Poisson process $X$   whose intensity may take on two possible values $\l_0$ or $\l_1$ (assume $\l_1>\l_0>0$ without loss of generality). 
 Given that the process $X$ is observed in real time, the problem is to detect the correct intensity as soon as possible and with minimal probabilities of incorrect terminal decisions.
  The above problem is usually called {\it sequential testing problems} for two simple hypotheses which has been solved in \cite{PesShir2000} if a  priori distribution on the true intensity is given (the probability of $\lambda=\lambda_0$ is $1-\pi$ and the probability of $\lambda=\lambda_1$ is $\pi$). 
  The key step of solving such a sequential testing problem lies in deriving an equivalent optimal stopping problem through a Bayesian formulation. Then the optimal stopping time obtained for the optimal stopping problem determines an optimal testing rule for the sequential testing problem subsequently.
Using the similar idea,  \cite{Daynik2006,Daynik2008} solve the sequential testing problems for compound Poisson process subsequently.\smallskip

Besides Poisson process, there also exists a well-established body of literature on the so-called Wiener sequential testing problem for which the observed process $X$ satisfies a stochastic differential equation.
For example,  to test the true drift(s) of  a Brownian motion or a stochastic differential equation with Gaussian noise, the testing problems for one-dimensional or multidimensional observed processes have been studied in \cite{Buo2016,
ErPZ2020,Gape2004,Gape2011,JP2018,JP2021}. It is worth mentioning that   Bayesian formulations are also adopted in those papers thanks to the priori distribution.\smallskip

In practice, it is more common to see that the priori distribution is not given. This paper is in part motivated by the following question: is it possible (within the Bayesian setting) to reformulate the statement of the sequential testing problem so that the solution to the optimal stopping problem does not depend on knowledge of the a priori distribution $\pi$?  To answer this question, one considers its minimax version by taking into consideration of all possible priori distributions.  The prior distribution of $\pi$ which yields the worst scenario among all possible a priori distributions is called the \textit{least favorable distribution} (see, for example, \cite{Lehmann} and references therein).\smallskip

The study on minimax testing problems could date back to \cite{Hu1965} and then has aroused great interests of mathematicians afterward. The readers are referred to the  review paper  \cite{Fa2021} and references therein for the recent developments if the observation is a discrete-time sequence and admits some general Markovian structure. Compared to the sequential testing problems with a priori distribution, the additional efforts for minimax sequential testing problems are devoted to identifying the least favorable distribution first. After a least favorable distribution is discovered, the optimal testing rule can be determined by a sequential testing problem directly. For the most recent developments in the discrete-time models,  one may refer to \cite{Fa2016} when the observations are i.i.d and to \cite{Fa2020} when the observations admit a general Markovian structure.   In the continuous-time model, when the observed process is a continuous-time process satisfying a stochastic differential equation, the minimax sequential testing problem for the drift has been solved in  \cite{ErMei2022}. 
In this paper, we aim to extend the minimax sequential testing theory for the intensity of a Poisson process. Similar to  \cite{Fa2016,Fa2020,ErMei2022}, our main effort is devoted to deriving the least favorable distributions.  To the best knowledge of the authors, the problem has not been considered in the literature.  \smallskip

We will follow the approaches adopted in  \cite{ErMei2022} while several new features have to be addressed in this paper.  First, the sample paths of  Poisson process are not continuous. Second, when working on the optimal stopping problem, we will see that the probabilistic regularity of the optimal stopping boundary, which is an essential step in \cite{ErMei2022}, fails in our case. Some new methods are needed here. Third, we will see that 1/2 and 1/2 is not the least favorable priori distribution even if the performance functional is symmetric. Such a result implies that the minimax sequential testing problem for Poisson process is not symmetric,  different from the linear minimax Wiener sequential testing problems in \cite{ErMei2022} for stochastic differential equations. Working on those new features, we will present an equivalent characterization of the least favorable distribution and present a sufficient condition for its existence.  Moreover, we will propose a numerical method for a least favorable distribution from its characterization. \smallskip

The rest of the paper is arranged as follows. We will present the explicit formulation of our minimax sequential testing problem in Section \ref{sec:for}. Then the optimal stopping theory of sequential testing problem is recalled in Section \ref{sec:stp} and some new results are introduced as a preparation of minimax sequential testing problem. In Section \ref{sec:mst}, we present our main results for minimax sequential testing problem: deriving an equivalent characterization for the least favorable distribution and concluding a sufficient condition for its existence. Then an illustrative example is presented in Section \ref{sec:ill} to explain why a different and easier formulation of the minimax problem is problematic. 
Then a numerical example is presented in Section \ref{sec:exp} to illustrate our main results. Finally, some concluding remarks are made in Section \ref{sec:con}.
	
	\section{Formulation}\label{sec:for}
Recall that the intensity $\l$ of Poisson process $X$ may take either $\l_1>\l_0>0$.	 The detection  problem is to test the following two simple hypotheses:
	$$H_0: \l=\l_0\text{ verses } H_1: \l=\l_1.$$
	We will formulate the sequential testing problem first and then the minimax sequential testing problem.\smallskip

1. We adopt the Bayesian formulation in \cite{PesShir2000} and assume the intensity $\lambda$  admits a priori distribution 
\bel{prior0}\PP(\l=\l_0)=1-\pi\text{ and }\PP(\l=\l_1)=\pi\eel
	for some  $\pi\in(0,1)$.
 We write $\psi=\pi/(1-\pi)$. One can easily see that $\pi$ and $\psi$ are essentially equivalent. Based on such a priori distribution, we define a new probability measure $\PP_\psi$ by
 \bel{ppvarphi}\PP_\psi(\cdot)=\frac {1}{1+\psi}\PP(\cdot|\l=\l_0)+\frac {\psi}{1+\psi}\PP(\cdot|\l=\l_1)\eel
Then 	\bel{prior}\PP_\psi(\l=\l_0)=\frac {1}{1+\psi}\text{ and }\PP_\psi(\l=\l_1)=\frac {\psi}{1+\psi}\eel
	for some  $\psi\in(0,\infty).$ 
Such a priori distribution only depends on a parameter $\psi$. Therefore we will call $\psi$ as the priori distribution of $\l$ if its distribution satisfies \eqref{prior}.\smallskip

The sequential testing problem aims at detecting the true parameter $\lambda$ based on the continuous observation of the state process. This yields that a decision to make consists of two perspectives:  a stopping time $\tau$ and the decided index $d=0$ or $1$ (i.e. 0 stands for accepting $H_0$ and 1 stands for accepting $H_1$).	To evaluate the performance of the decision couple $(\tau,d)$, we define the performance functional  by a linear combination of {\it the expected waiting time} and {\it the expected probabilities of false detection}:
\bel{per}J(\psi;\tau,d)=\EE_{\psi}\big[\tau+aI(d=0; \l=\l_1)+bI(d=1; \l=\l_0)\big]\eel
	for some given $a,b>0$. The two constants $a,b$ stand for the costs for two possible types of wrong detection. If $a=b$, the cost functional is called {\it symmetric} as the costs of false detection are the same.
  The sequential testing problem is  to find an optimal $(\tau^*,d^*)$ such that the performance functional $J(\cdot)$ is minimized.\smallskip

  2. The above sequential testing problem has been solved in \cite{PesShir2000}. The main approach is to derive an equivalent optimal stopping problem through a measure-change method. Let us summarize the key steps as follows. 
  \smallskip
  
  Write $\dbS$ by the set of all stopping time $X$ and $\cF_t^X$ by the natural filtration of $X$.   {\it The likelihood ratio process } $L_t$,  {\it  the posterior probability process} $\Pi_t$, and {\it the posterior probability ratio  process} $\varPsi_t$ can be written by 
 $$L_t=\frac{d\PP(\cdot|\l=\l_1)}{d\PP(\cdot|\l=\l_0)}\Big|_{\cF_t^X}, ~\Pi_t:=\PP(\lambda=\lambda_1|\cF_t^X)\text{ and } \varPsi_t=\frac{\Pi_t}{1-\Pi_t}.$$
From (2.14) and (2.16) in \cite{PesShir2000}, we see that the posterior probability ratio  process $\varPsi$ satisfies
\begin{align}\label{eqphi} &\varPsi_t=\psi L_t=\psi\exp\left\{(X_t-X_0)\log\(\frac{\l_1}{\l_0}\)-(\l_1-\l_0)t\right\},\nonumber\\& d\varPsi_t=\(\frac{\l_1}{\l_0}-1\)\varPsi_{t^-}d(X_t-\l_0t).\end{align}
Similar to  (2.9) in \cite{ErPZ2020}, we are able to derive a new form of the $\bar J$ by changing the measure from $\EE_\psi$ to $\EE_0$. 
First we note that \begin{align*}
&J(\psi;\tau,d)=\EE_{\psi}\big[\tau+aI(d=0; \l=\l_1)+bI(d=1; \l=\l_0)\big]\\
&\q=\EE_\psi[\tau+a\Pi_t I(d=0)+b(1-\Pi_t)I(d=1)]\\
&\q\geq \EE_\psi\(\tau+[a\Pi_\t] \wedge [b(1-\Pi_\t)]\)=:\bar J(\psi;\tau)
\end{align*}
The equality holds if and only 
$d=0$ if $a\Pi_\t< b\Pi_\t$, $d=1$ if $a\Pi_\t> b\pi_\t$
and $d=0 $ or $1$ if $a\Pi_\t=b\Pi_\t$.
Now it follows that 
\begin{align}\label{cost0000}
&\bar J(\psi;\tau)= \EE_\psi\(\tau+[a\Pi_\t] \wedge [b(1-\Pi_\t)]\)\nonumber\\
&=\frac{1}{1+\psi}\EE_0\(\int_0^\tau ds+[a\Pi_\t] \wedge [b(1-\Pi_\t)]\)+\frac{\psi}{1+\psi}\EE_1\(\int_0^\tau ds+[a\Pi_\t] \wedge [b(1-\Pi_\t)]\)\nonumber\\
&\q=\frac{1}{1+\psi}\EE_0\(\int_0^\tau 1 ds+[a\Pi_\t] \wedge [b(1-\Pi_\t)]\)\nonumber\\
&\qq\qq+\frac{\psi}{1+\psi}\EE_0\(\int_0^\tau L_s ds+L_\tau\times [a\Pi_\t] \wedge [b(1-\Pi_\t)]\)\nonumber\\
&=\frac{1}{1+\psi}\EE_0\(\int_0^\tau (1+\varPsi_s) ds+(1+\varPsi_\tau)\times [a\Pi_\t] \wedge [b(1-\Pi_\t)]\)
\nonumber\\
&\q=\frac{1}{1+\psi}\EE_0\Big[\int_0^\tau(1+\varPsi_t)dt+(a\varPsi_\t)\wedge b\Big].\end{align}
We note that $X_t$ is a Poisson process with intensity $\lambda_0$ under $\PP_0$ and subsequently $\varPsi$ is a strong Markov process under $\PP_0$. Combining \eqref{eqphi} with \eqref{cost0000}, we obtained the equivalent optimal stopping problem for the sequential testing problem if the priori distribution is given:\smallskip

\noindent{\bf Problem (ST):}    Find a $\tau^*\in \dbS$ such that 
\bel{seqt}\bar J(\psi;\tau^*)=\inf_{\tau\in\dbS}\bar J(\psi;\tau).\eel

3. For {\bf Problem (ST)}, it is proved  in \cite{PesShir2000} that there exists two positive numbers $\a^*$ and $\beta^*$ such that the optimal stopping time is $$\tau^*=\inf\{t\geq 0:\varPsi_t\notin(\a^*,\beta^*)\}$$
and the optimal decision $d^*$ at $\tau^*$ satisfies
$$d^*_{\t^*}=\left\{\ba{ll}1,\qq\qq\q\text{ if } \varPsi_{\tau^*}>b/a;\\
0,\qq\qq\q\text{ if } \varPsi_{\tau^*}<b/a;\\
\text{either 0 or 1},\text{ if } \varPsi_{\tau^*}=b/a.\ea\right.$$
It is worth mentioning that $\bar J(\psi;\tau)=J(\psi;\tau,d_\tau^*)$.\smallskip

4. With the above preparation in hand, we turn to the formulation of the minimax sequential testing problem. The objective of the minimax formulation is to minimize the performance functional $\bar J(\psi;\cdot)$ in the \textit{worst case scenario} of all prior distributions $\psi$.  This leads to the following optimal stopping problem. \medskip
	
	\noindent{\bf Problem MinimaxST:}   Find a $(\psi^*;\tau^*)\in (0,\infty)\times\dbS$ such that
	\bel{minmaxopt}\bar J(\psi^*;\tau^*)=\inf_{\tau\in\dbS}\sup_{\psi> 0}\bar J(\psi;\tau).\eel

5. The optimization problem in \eqref{minmaxopt} belongs to the class of minimax optimization problems.	We aim to find a  saddle-point $(\psi^*,\tau^*)$ for {\bf Problem MinimaxST}, i.e. 
\bel{saddle}\bar J(\psi;\tau^*)\leq  \bar J(\psi^*;\tau^*)\leq\bar J(\psi^*;\tau)\text{ for any }(\psi,\tau)\in (0,\infty)\times \dbS.\eel
	The saddle-point property consists of  two perspectives: \smallskip
	
	(1) $\tau^*$ is the optimal stopping time for {\bf Problem ST} if $\psi=\psi^*$.  This can be concluded from the sequential testing problem directly. We will still need to conclude some auxiliary results for the future minimax sequential testing problem.\smallskip

	(2) Given  $\tau^*\in \dbS$, $\psi=\psi^*$ is the {\it least favorable distribution} in $\bar J(\psi;\tau^*)$, i.e. $$\bar J(\psi;\tau^*)\leq \bar J(\psi^*;\tau^*)\text{ for any }\psi\in(0,\infty).$$ 
	This is our main focus in the paper. The idea is to take the derivative of $\bar J(\cd;\tau^*)$ and investigate its monotonicity on $[0,\infty)$. We will see this later.\smallskip

 This paper's key result for the minimax sequential testing problem for Poisson processes, given by Theorem \ref{mainthh}, is that there exists a $\varphi_0$ which is the least favorable distribution for the stopping time $\tau^*(\varphi_0)$. We want to emphasize that $\varphi_0=1$ is not necessarily the least favorable distribution even if the performance functional is symmetric (i.e. $a=b$). This is to say the intuitive $1/2$-$1/2$ rule for symmetric costs is not true in general. This will be justified by a numerical example in Section \ref{sec:exp}. 

	\section{Sequential testing and some auxiliary results}\label{sec:stp}
In this section, we will recall the optimal stopping theory for sequential testing problems for the intensity of  Poisson processes from \cite{PesShir2000}. Additionally, we will present some auxiliary results which are useful for our minimax sequential testing problem in the next section.\smallskip

The following theorem characterizes the optimal stopping boundary which is a consequence of Theorem 2.1 in \cite{PesShir2000}. 
\begin{theorem}\label{STthm} Suppose the priori distribution $\psi$ is given.\smallskip

{\rm (1)} If \bel{l1l00}\l_1-\l_0\leq\frac1a+\frac1b,\eel
the stopping time $\tau^*=0$ is optimal for the sequential testing problem in \eqref{seqt}.\smallskip

{\rm (2) }If \bel{l1l0}\l_1-\l_0>\frac1a+\frac1b,\eel
then the stopping time 
$$\tau^*=\inf\Big\{t\geq 0:\psi L_t\notin (\alpha^*,\beta^*)\Big\}$$
is optimal for the sequential testing problem in \eqref{seqt} for some unique $\alpha^*, \beta^*$ satisfying $0<\a^*<b/a<\b^*$.
\end{theorem}

In the above theorem, the closed forms of $\alpha^*$ and $\beta^*$ may not be obtainable, while they can be found through some numerical algorithms (see \cite{PesShir2000,Daynik2006}).  Therefore we will treat them as given constants in the sequel.  Because $\varPsi_t$ is strictly decreasing between two jumps,   the optimal lower stopping boundary $\a^*$ is  probabilistic regular while the optimal upper stopping boundary $\beta^*$ is not probabilistic regular (see Section 4.2 in \cite{KA1998}), i.e. for $\tau=\inf\{t\geq0: L_t>\beta^*\}$, it follows that
$$\PP_{\beta^*}(\tau>0)=1\neq 0.$$ The probabilistic regularity is a critical property used in \cite{ErMei2022} and thus we need to provide a different approach here.\smallskip

Now  define  $$\tau^*(\varphi_0):=\inf\Big\{t\geq 0:\varphi_0 L_t\notin (\alpha^*,\beta^*)\big\}.$$ Note that $L$ is decreasing between  jumps and $L$ only jumps to a larger value at each jump time.
At $t=\tau^*(\varphi_0)$, we have either $L_t=\alpha^*/\varphi_0$ or $L_t\geq \beta^*/\varphi_0$. 
We proceed with the following proposition.
\begin{proposition}\label{conttau} {\rm (1)} For any $\varphi_0 \in (\alpha^*,\beta^*) $,
$$\tau^*(\varphi)\rightarrow \tau^*(\varphi_0),~~\PP_0\text{ - a.s. as } \varphi\rightarrow\varphi_0.$$

{\rm (2)} It follows that 
$$\tau^*(\varphi)\rightarrow 0,~~\PP_0\text{ - a.s. as }  \varphi\downarrow \alpha^*.$$

{\rm (3)} It follows that $$\lim_{\varphi\uparrow\beta^*}\tau^*(\varphi)=\inf\Big\{t\geq 0: L_t\notin (\alpha^*/\beta^*,1]\Big\}=:\th(\beta^*)~\PP_0\text{-a.s.} .$$
\end{proposition}
\begin{proof}
(1) Note that $$\tau^*(\varphi):=\inf\Big\{t\geq 0: L_t\notin (\alpha^*/\varphi,\beta^*/\varphi)\Big\}.
$$
 Because $L$ is  left-continuous on stopping times a.s. and is  right-continuous,   we have 
\begin{align*}&\lim_{\varphi\downarrow\varphi_0}\tau^*(\varphi)=\inf\Big\{t\geq 0: L_t\notin [\alpha^*/\varphi_0,\beta^*/\varphi_0)\Big\},\\
&\lim_{\varphi\uparrow\varphi_0}\tau^*(\varphi)=\inf\Big\{t\geq 0: L_t\notin (\alpha^*/\varphi_0,\beta^*/\varphi_0]\Big\}=:\th(\varphi_0).\end{align*}
Because  $\a^*/\varphi_0$ is a probabilistic boundary for  $[\alpha^*/\varphi_0,\beta^*/\varphi_0)$,  it follows that  $\lim_{\varphi\downarrow\varphi_0}\tau^*(\varphi)=\tau^*(\varphi_0)$.
Note that $\tau(\varphi_0)\leq \th(\varphi_0)$, it suffices to show that 
\bel{P0}\PP_0(L_{ \tau(\varphi_0)}=\beta^*/\varphi_0)=0.\eel

Since $L_t$ is strictly decreasing between each jump and there are at most countable many jumps for each possible path of $L$,   $\{L_t=\ell\}$ happens only for countable many $t$ for any $\ell\geq 0$. Since $\beta^*/\varphi_0>1$, before $\th(\varphi_0)$, the only possibility for $L_t$ attaining $\beta^*/\varphi_0$ is that a jump happens exactly at the time when $L_{t}$ hits $\beta^*\l_0/(\varphi_0\l_1)$. As $L_t$ hits $\beta^*\l_0/(\varphi_0\l_1)$ at most countable many times on the whole time horizon $[0,\infty)$ and the probability that a jump happened at those stopping times are $0$, this yields  \eqref{P0}. Therefore we have proved that 
$$\lim_{\varphi\rightarrow\varphi_0}\tau^*(\varphi)=\tau^*(\varphi_0)~\PP_0\text{-a.s.} .$$
It is evident that (2) and (3)  are special cases of (1). The proof is complete.
\end{proof}

\begin{remark}\rm 
We remark that the above results do not hold for $\varphi_0=\beta^*$ because $\tau^*(\beta^*)=0$ a.s.  and 
$\lim_{\varphi\uparrow\beta^*}\tau^*(\varphi)>0$ a.s. (we need to wait for the first jump).\end{remark}

In order to work on the minimax sequential testing problem in the next section,  we proceed with the following two differential equations: 
   \bel{ode1}\left\{\ba{ll}\ns\ds (\l_0-\l_1)\varphi f_0'(\varphi)+\l_0[f_0(\l_1\varphi/\l_0)-f_0(\varphi)]=0 \text{ in }(\alpha^*/\beta^*,1),\\
   \ns\ds f_0(1)=1; \q f_0(\varphi)=0,\text{ for }\varphi> 1.\ea\right.\eel
   and 
   \bel{ode2}\left\{\ba{ll}\ns\ds (\l_0-\l_1)\varphi f_1'(\varphi)+\l_0[f_1(\l_1\varphi/\l_0)-f_1(\varphi)]+(\varphi-1)=0 \text{ in }(\alpha^*/\beta^*,1),\\
   \ns\ds f_1(1)=0; \q f_1(\varphi)=-b,\text{ for }\varphi> 1.\ea\right.\eel
Note that  $\l_1>\l_0>0$, it can be easily seen that the above two differential equations admit unique continuous solutions $f_0$ and $f_1$ on $[\a^*/\beta^*, 1]$. We  can also conclude the monotonicity of $f_0$ and $f_1$ as follows.
\begin{lemma} The solutions $f_0$ and $f_1$ to \eqref{ode1} and \eqref{ode2} respectively are strictly decreasing on  $[\a^*/\beta^*, 1]$.
\end{lemma}
\begin{proof} 

Note that  $\l_1>\l_0>0$. Suppose that  $\varphi_0$ is the largest $\varphi $ in $(\a^*/\b^*,1)$ such that
$f_0'(\varphi_0)=0$. This is to say that $f_0'(\varphi)<0$ for any $\varphi\in (\varphi_0,1)$ which says that $f$ is strictly decreasing on $(\varphi_0,1]$. In this case, $f_0'(\varphi_0)=0$ leads to
$$f_0(\l_1\varphi_0/\l_0)-f_0(\varphi_0)=0,$$
which is a contradiction. Therefore  $f'(\varphi)<0$ for all $\varphi \in(\a^*/\b^*,1)$. The proof for $f_1$ is the same and thus omitted here.
\end{proof}

With above lemma, we  define 
\bel{gamma*}\gamma^*:=\frac{ a\alpha^*/\beta^*-f_1(\alpha^*/\beta^*)}{f_0(\alpha^*/\beta^*)}.\eel
Note that $\gamma^*$ is well-defined because $f_0(\alpha^*/\beta^*)\geq f_0(1)=1.$ We proceed with the following lemma.
\begin{lemma}\label{thbeta}
Recall $$\th(\beta^*)=\inf\Big\{t\geq 0: L_t\notin (\alpha^*/\beta^*,1]\Big\}\text{ with }L_0=1.$$
Then it follows that
$$\gamma^*=\(\frac{a\alpha^*}{\beta^*}+b\)\PP_0(L_{\th(\beta^*)}=\alpha^*/\beta^*)-b+\EE_0\int_0^{\th(\beta^*)}(L_t-1)dt.$$

\end{lemma}

\begin{proof} The proof is based on applying It\^o's formula to $$t\mapsto \gamma^*f_0(L_t)+f_1(L_t).$$ By Dynkin's formula, noting that $\PP_0(L_{\th(\beta^*)}=1)=0$ (see \eqref{P0}), we have 
$$\ba{ll}\ds\EE_0\int_0^{\th(\beta^*)}(L_t-1)dt=\gamma^*-\EE_0\(\gamma^*f_0(L_{\th(\beta^*)})+f_1(L_{\th(\beta^*)})\)\\
\ns\ds= \gamma^*- [\gamma^* f_0(\alpha^*/\beta^*)+f_1(\alpha^*/\beta^*)]\PP_0(L_{\th(\beta^*)}=\alpha^*/\beta^*)+b(1-\PP_0(L_{\th(\beta^*)}=\alpha^*/\beta^*))\\
\ns\ds=\gamma^*-(a \alpha^*/\beta^*+b)\PP_0(L_{\th(\beta^*)}=\alpha^*/\beta^*)+b.\ea$$
The proof is complete.
\end{proof}
    
   \section{Minimax sequential testing}\label{sec:mst}

   In this section, we will solve the minimax sequential testing problem. Our goal is to characterize the least favorable distribution and prove its existence. \smallskip
   
   As proved in the previous section,  $\tau^*(\varphi_0)$ is the optimal stopping time if the priori distribution is $\varphi_0$. We aim to characterize the least favorable distribution $\varphi_0$ in the sense that 
   $$\bar J(\psi;\tau^*(\varphi_0))\leq \bar J(\varphi_0;\tau^*(\varphi_0))\text{ for any }\psi>0.$$
   From Theorem \ref{STthm},  we have two different cases.\smallskip
   
   (1) When \eqref{l1l00} holds, $\tau^*(\varphi_0)=0$ is the optimal decision for any $\varphi_0$. In this case, we have
   \bel{barj1}\bar J(\psi;0)=\frac{b\wedge (a\psi)}{1+\psi}.\eel
   Therefore the maximum point $\varphi_0=b/a$ of $\bar J(\cdot;0)$ is the unique least favorable distribution. \smallskip

   (2) Suppose \eqref{l1l0} holds. In this case, we have
   $$\bar J(\psi;\tau^*(\varphi_0))=\frac1{1+\psi}\EE_0\Big[\int_0^{\tau^*(\varphi_0)}(1+\psi L_t)dt+b \wedge(a\psi L_{\tau^*(\varphi_0)})\Big].$$

First, we claim that the least favorable distribution $\varphi_0$ must fall in $(\a^*,\beta^*)$. Suppose $\varphi_0\notin (\a^*,\beta^*)$. Then we have $\tau^*(\varphi_0)=0$ and $\bar J(\psi;0)$ takes the same  form as \eqref{barj1} whose maximum point is $b/a\in(\a^*,\beta^*)$. Therefore $\varphi_0$ can not be a least favorable distribution.\smallskip

 Then we only need to consider $\varphi_0\in (\a^*,\beta^*)$.
 Our idea is to find the (left and right) derivatives of  $\bar J(\cdot;\tau(\varphi_0))$. Note that 
$$\ba{ll}\ad\lim_{\delta\rightarrow 0^+}\frac1\delta\Big[\EE_{0}\(b\wedge (a(\psi+\delta) L_{\tau^*(\varphi_0)})\)-\EE_{0}\(b\wedge (a\psi L_{\tau^*(\varphi_0)})\)\Big]\\
\ns\ad=\lim_{\delta\rightarrow 0^+}\frac 1 \delta\Big[a\delta\EE_0 \(L_{\tau^*(\varphi_0)}I[ L_{\tau^*(\varphi_0)}\leq \frac b{a(\psi+\delta)}]\)\\
\ns\ad\qq\qq+\EE_0 \((b-a\psi L_{\tau^*(\varphi_0)})I[ \frac b{a(\psi+\delta)}< L_{\tau^*(\varphi_0)}<\frac b{a\psi}]\)\Big]\\
\ns\ad=a\EE_0 \(L_{\tau^*(\varphi_0)}I[ L_{\tau^*(\varphi_0)}< \frac b{a\psi}]\).\ea$$
Taking the right derivative of  $\bar J(\cd;\tau(\varphi_0))$, we have 
$$\ba{ll}\ad \frac{\partial^+ \bar J}{\partial \psi }(\psi;\tau^*(\varphi_0))\\
\ns\ad=\frac1{(1+\psi)^2}\Big[a(1+\psi)\EE_0\(L_{\tau^*(\varphi_0)}I(\psi L_{\tau^*(\varphi_0)}<b/a )\)\\
\ns\ad\qq\qq-\EE_0 \(b\wedge (a\psi L_{\tau^*(\varphi_0))}\)-\EE_0\int_0^{\tau^*(\varphi_0)}(1-L_t) dt\Big]\\
\ns\ad=\frac1{(1+\psi)^2}\Big[a\EE_0\(L_{\tau^*(\varphi_0)}I[\psi L_{\tau^*(\varphi_0)}<b/a]\) \\
\ns\ad\qq\qq-b\PP_0 (\psi L_{\tau^*(\varphi_0)}\geq b/a)-\EE_0\int_0^{\tau^*(\varphi_0)}(1-L_t)dt\Big] \ea$$
Note that 
\begin{equation*}\ba{ll}\ad\lim_{\delta\rightarrow 0^+}\frac1\delta\Big[\EE_0 [b\wedge (a\psi L_{\tau^*(\varphi_0)})]-\EE_0 [b\wedge (a(\psi-\delta) L_{\tau^*(\varphi_0)})]\Big]\\
\ns\ad=\lim_{\delta\rightarrow 0^+}\frac1\delta\Big[a\delta\EE_0 \(L_{\tau^*(\varphi_0)}I[ L_{\tau^*(\varphi_0)}\leq \frac b{a\psi}]\)\\
\ns\ad\qq\qq+\EE_0 \((b-a\psi L_{\tau^*(\varphi_0)})I[ \frac b{a\psi}< L_{\tau^*(\varphi_0)}<\frac b{a(\psi-\delta)}]\)\Big]\\
\ns\ad=a\EE_0 \(L_{\tau^*(\varphi_0)}I[ L_{\tau^*(\varphi_0)}\leq  \frac b{a\psi}]\).\ea\end{equation*}
Similarly, taking  the left derivative with respect to $\psi$, we have
\begin{equation*}\ba{ll}\ad \frac{\partial^- \bar J}{\partial \psi }(\psi;\tau^*(\varphi_0))=\frac1{(1+\psi)^2}\Big[a\EE_0\(L_{\tau^*(\varphi_0)}I(\psi L_{\tau^*(\varphi_0)}\leq b/a)\)\\
\ns\ad\qq\qq\qq-b\PP_0 (\psi L_{\tau^*(\varphi_0)}> b/a)-\EE_0\int_0^{\tau^*(\varphi_0)}(1-L_t)dt\Big].\ea\end{equation*} 

 We notice that that the $(1+\psi)^2\frac{\partial^\pm \bar J}{\partial \psi }(\psi;\tau^*(\varphi_0))$ are decreasing on $(\a^*,\b^*)$ and for $\varphi_0\in (a^*,\b^*)$,  \bel{conh}\PP_0(\varphi_0 L_{\tau^*(\varphi_0)}=b/a)=0\text{ and 
} \varphi_0L_{\tau^*(\varphi_0)}=\alpha^* \text{ on } \{\varphi_0L_{\tau^*(\varphi_0)}<b/a \}.\eel Therefore, $\bar J(\cdot;\tau^*(\varphi_0))$ takes a maximum at $\varphi_0$ if and only if 
   \begin{equation*} h(\varphi_0):=\(\frac{a \alpha^*}{\varphi_0}+b\) \PP_0(\varphi_0 L_{\tau^*(\varphi_0)}<b/a) -b+\EE_0\int_0^{\tau^*(\varphi_0)}(L_t-1)dt=0.\end{equation*}
   By Proposition \ref{conttau} and \eqref{conh}, invoking dominant convergence, we know  $h(\cdot)$ is continuous with respect to $\varphi_0\in [\alpha^*,\beta^*)$ with $$ h(\alpha^*)=a>0.$$
Note that $\tau^*(\varphi_0)\rightarrow\th(\b^*)$ as $\varphi_0\uparrow \b^*$.  By the definition of $\gamma^*$ and \eqref{conh}, it also follows that \bel{hb*}h({\beta^*}^-)=\(\frac{a\alpha^*}{\beta^*}+b\)\PP_0(\b^*L_{\th(\b^*)}<b/a)-b+\EE_0\int_0^{\th(\b^*)}(L_t-1)dt=\gamma^*.\eel
By the mean-value theorem, all the above calculations yield the main results for the minimax sequential testing problem as follows.
   \begin{theorem} \label{mainthh}
{\rm (1)} If $\l_1-\l_0\leq1/a+1/b,$
the stopping time $\tau^*=0$ is the  optimal stopping time for the minimax sequential testing problem in \eqref{minmaxopt} where $\varphi_0=b/a$ is the least favorable distribution.\smallskip

{\rm (2) }If $\l_1-\l_0>1/a+1/b,$ the stopping time $\tau^*(\varphi_0)$ is the optimal stopping time for the minimax sequential testing problem in \eqref{minmaxopt} if and only if
   $h(\varphi_0)=0.$
   Moreover, if $\gamma^*<0$, there exists at least one least favorable distribution $\varphi_0\in (\alpha^*,\beta^*).$ 
   \end{theorem}

If the saddle point exists, it is well-known that the strong duality for the minimax problem holds, i.e. $$\inf_{\tau\in\dbS}\sup_{\psi> 0}\bar J(\psi;\tau)=\sup_{\psi> 0}\inf_{\tau\in\dbS}\bar J(\psi;\tau)=\sup_{\psi> 0}V(\psi)$$
where $V(\cdot)$ is the value function of the sequential testing problem. This says that the least favorable distribution $\psi^*$ is the maximum point of $V(\cdot)$. Observing from above, it immediately holds the following corollary.
\begin{corollary}
If $\l_1-\l_0\leq1/a+1/b$ or if $\l_1-\l_0>1/a+1/b$ and $\gamma^*<0$, then the maximum point of $V(\psi)=\inf_{\tau\in\dbS}\bar J(\psi;\tau)$ is the least favorable distribution.
\end{corollary}
Finally, let us conclude the section with a numerical algorithm to find the least favorable distribution $\varphi_0$ if $\l_1-\l_0>1/a+1/b$. (1) We find the optimal stopping boundary $\a^*$ and $\b^*$ numerically. (2) We find $\gamma^*$  through \eqref{gamma*}. (3) If $\gamma^*<0$, we can proceed the bisection search for the zero-point of $h$ on $(\a^*,\b^*)$  or find the maximum of the value function $V$.

\section{Illustrative Example}\label{sec:ill}

In this section, we present an example to illustrate why a different and easier formulation for our  {\bf Problem MinimaxST} is problematic.

We note that the main goal of this paper is to find a saddle point $(\varphi^*,\tau^*)$  for {\bf Problem MinimaxST}, where $\tau^*$ is a stopping time. It is well-known that a stopping time can be defined by some appropriate stopping policies on the paths of the underlying processes. Some readers may suggest to reformulate our minimax problem as follows.

Write $\vec \varPsi^\psi$ by the path of $\varPsi_t$ with initial $\psi$. Define $\Xi$ by the set of all possible stopping policies.
\medskip

\noindent{\bf Problem Minimax-2:}
$$\inf_{\xi\in\Xi}\sup_{\psi\geq 0}\bar J(\psi;\xi(\vec \varPsi^\psi)). $$

We note that in the above problem, the stopping time $\tau$ is replaced by applying the stopping policy on the path of $\varPsi_t$.

(1) {\bf Problem Minimax-2} is trivial because $\xi_*$, defined by the first time enters the stopping region  $D$, is optimal for all $\psi$. Then it follows that
$$\inf_{\xi\in\Xi}\sup_{\psi\geq 0}\bar J(\psi;\xi(\vec \varPsi^\psi))=\sup_{\psi\geq 0}\bar J(\psi;\xi_*(\vec \varPsi^\psi))=\sup_{\psi\geq 0}\inf_{\xi\in\Xi}\bar J(\psi;\xi(\vec \varPsi^\psi)). $$
Therefore it can be claimed that $(\psi_*,\xi_*)\in [0,\infty)\times \Xi$ is the saddle point for  {\bf Problem Minimax-2}, where  $\psi_*$ is the maximum point of $\psi\mapsto J(\psi;\xi_*(\vec \varPsi^\psi))$. Therefore   {\bf Problem Minimax-2} is much easier to solve than  {\bf Problem MinimaxST} proposed in the paper. Now it is very natural to ask why it is necessary to consider {\bf Problem MinimaxST} in the paper. In fact, solving {\bf Problem Minimax-2} only will be problematic in some simple cases.

(2) Now let us present a simple example to explain why solving {\bf Problem Minimax-2} is problematic.

Let $X_t=x-t$ under $\dbP_x.$ Note that the path of $X_t$ can be determined by the initial value only, we write $\vec X=\vec x$ when $X_0=x$. Because the filtration is trivial, the set of all stopping times is $[0,\infty)$.  Adopting the ideas of {\bf Problem MinimaxST} and {\bf Problem Minimax-2}, we have the following two problems coorespondingly:\medskip

\noindent {\bf Problem Minimax-3:}
$$\inf_{\tau\in[0,\infty)}\sup_{x\geq0}\dbE_x|X_\tau|=\inf_{\tau\in[0,\infty)}\sup_{x\geq0}|x-\tau|.$$

\noindent {\bf Problem Minimax-4:}
$$\inf_{\xi\in\Xi}\sup_{x\geq0}\dbE_x|X_{\xi(\vec X)}|=\inf_{\xi\in\Xi}\sup_{x\geq0}|x-\xi(\vec x)|.$$

It can be easily seen that {\bf Problem Minimax-3} admits no saddle point in $[0,\infty)\times [0,\infty)$. While $(x,\xi_*)\in [0,\infty)\times \Xi$ is always a saddle point of {\bf Problem Minimax-4} for all $x$ where $\xi_*$ is defined by the first time to hit $0$. {\bf Problem Minimax-4} says that any initial value of $x$ can be a good choice which is against the goal of deriving a robust decision when proposing a minimax optimization problem. Instead, the result from {\bf Problem Minimax-3} says that there is no robust choice in such a problem.

(3) In fact, when {\bf Problem MinimaxST} admits a saddle point $(\varphi_*,\tau_*)$, then $(\varphi_*,\xi_*)$ is a saddle point {\bf Problem Minimax-2}. Even though the paper is devoted to proving  {\bf Problem MinimaxST} admits a saddle point, solving the easier version {\bf Problem Minimax-2} is still problematic in our paper.

\section{Numerical Example}\label{sec:exp}
In this section, we present a concrete example to illustrate our theory. The example is also used in  \cite{PesShir2000} (see Figure 2 there). Suppose
$$\l_0=1,~\l_1=5,~a=b=2.$$
In such case $\l_1-\l_0>1/a+1/b,$ and we have
$$L_t=\exp\Big\{(X_t-X_0)\log 5-4t\Big\}.$$
The observed process $X_t$ has intensity $\l_0=1$ under $\PP_0$. The numerical method yields that
$$\a^*\thickapprox0.297,~~\b^*\thickapprox2.390,~~\frac{\a^*}{\b^*}\thickapprox0.124.$$
Through a numerical method, we find that 
$\gamma^*= -0.788<0.$
This concludes the existence of a least favorable distribution $\varphi_0$. A  bisection search yields the least favorable distribution $\varphi_0=0.977$. Therefore, the optimal stopping time for {\bf MinimaxST} is
$$\tau^*=\inf\Big\{t\geq 0: L_t\notin (0.304,2.440)\Big\}.$$
Moreover, we find that $h(1)=-0.03\neq 0.$ Therefore the minimax sequential testing problem is not symmetric even if $a=b$.
   
\section{Concluding remarks}\label{sec:con}
In this paper, we solved the minimax sequential testing problem for the intensity of a Poisson process. We find an equivalent characterization of the least favorable distribution and present a sufficient condition for its existence. Moreover, a numerical algorithm is proposed for the least favorable distribution. For further study on minimax sequential testing problems of continuous-time models, we would like to consider the minimax sequential problems for more general processes such as compound Poisson processes and L\'evy processes.

\newpage 

{\bf Statements and Declarations}

The first-named author is partially supported by Simons Foundation's Travel Support for Mathematicians Program (No. 00002835).

The authors have no relevant financial or non-financial interests to disclose.

\begin{thebibliography}{}

\bibitem{Buo2016} Buonaguidi, B., \& Muliere, P. (2016). Bayesian sequential testing for Lévy processes with diffusion and jump components. {\it Stochastics}, 88(7), 1099-1113.

\bibitem{Daynik2006} Dayanik, S., \& Sezer, S. O. (2006). Sequential testing of simple hypotheses about compound Poisson processes. {\it Stochastic processes and their applications}, 116(12), 1892-1919.
\bibitem{Daynik2008}Dayanik, S., Poor, H. V., \& Sezer, S. O. (2008). Sequential multi-hypothesis testing for compound Poisson processes. {\it Stochastics: An International Journal of Probability and Stochastic Processes}, 80(1), 19-50.
     
		
		\bibitem{ErPZ2020}Ernst, P. A., Peskir, G., \& Zhou, Q. (2020). Optimal real-time detection of a drifting Brownian coordinate. {\it The Annals of Applied Probability}, 30(3), 1032--1065.
		
	\bibitem{ErMei2022}  Ernst, P. A., \& Mei, H. (2025). The Minimax Wiener Sequential Testing Problem. {\it SIAM Journal on Control and Optimization}, 63(1), 206-226.
		\bibitem{Fa2016} Fauss, M., \& Zoubir, A. M. (2016). Old bands, new tracks—Revisiting the band model for robust hypothesis testing. {\it IEEE Transactions on Signal Processing}, 64(22), 5875-5886.
		\bibitem{Fa2020}	Fauss, M., Zoubir, A. M., \& Poor, H. V. (2020). Minimax optimal sequential hypothesis tests for Markov processes. {\it  The Annals of Statistics}, 48(5), 2599-2621.
		\bibitem{Fa2021}	Fauss, M., Zoubir, A. M., \& Poor, H. V. (2021). Minimax robust detection: Classic results and recent advances. 
 {\it  IEEE Transactions on Signal Processing}, 69, 2252-2283.
		
		
		
	
		
		\bibitem{Gape2004}	Gapeev, P. V., \& Peskir, G. (2004). The Wiener sequential testing problem with finite horizon. {\it  Stochastics and stochastic reports}, 76(1), 59--75.
		
		
		\bibitem{Gape2011}Gapeev, P. V., \& Shiryaev, A. N. (2011). On the sequential testing problem for some diffusion processes. {\it Stochastics: An International Journal of Probability and Stochastic Processes}, 83(4-6), 519--535.
		
	
		
		
		
		\bibitem{Hu1965} Huber, P. J. (1965). A robust version of the probability ratio test. {\it The Annals of Mathematical Statistics}, 1753-1758.
		
	
		
		
		\bibitem{JP2018}Johnson, P. \& Peskir, G. (2018). Sequential testing problems for Bessel processes.  {\it Transactions of the American Mathematical Society}, 370(3), 2085--2113.
		
		\bibitem{JP2021}	Johnson, P., Pedersen, J. L., Peskir, G., \& Zucca, C. (2022). Detecting the presence of a random drift in Brownian motion. {\it Stochastic Processes and their Applications}, 150, 1068--1090.
		
		\bibitem{KA1998}Karatzas, I., \& Shreve, S. E. (1998). {\it Brownian motion. In Brownian Motion and Stochastic Calculus} (pp. 47-127). Springer, New York, NY.
		
	\bibitem{Lehmann} Lehmann, E.L. (1959)
		{\it Testing Statistical Hypotheses}. John Wiley \& Sons. 


   \bibitem{PesShir2000}	Peskir, G., \& Shiryaev, A. N. (2000). Sequential testing problems for Poisson processes. {\it Annals of Statistics}, 837-859.



\bibitem{Wald1947} Wald, A. {\it Sequential Analysis.} Wiley, Hoboken, NJ, USA, 1947.
  
		
	
    	
    
    
    \end{thebibliography}
\end{document}